\documentstyle[11pt,twoside]{article}
\newtheorem{theorem}{\noindent\bf Theorem}[section]
\newtheorem{proposition}[theorem]{\noindent\bf Proposition}
\newtheorem{lemma}[theorem]{\noindent\bf Lemma}

\newtheorem{remark}[theorem]{\noindent\bf Remark}

\newtheorem{definition}[theorem]{\noindent\bf Definition}

\newcommand{\g}{{\cal G}}

\newcommand{\m}{\textsf{m}}
\newcommand{\n}{\textsf{n}}

\input{amssym}
\begin{document}
\date{}
\title{{\Large\bf Introverted subspaces of the duals of measure algebras }}

\author{{\normalsize\sc H. Javanshiri and R. Nasr-Isfahani}}

\maketitle
\normalsize

\begin{abstract}
Let $\g$ be a locally compact group. In continuation of our
studies on the first and second duals of measure algebras by the
use of the theory of generalised functions, here we study the
C$^*$-subalgebra $GL_0(\g)$ of $GL(\g)$ as an introverted subspace
of $M(\g)^*$. In the case where $\g$ is non-compact we show that
any topological left invariant mean on $GL(\g)$ lies in
$GL_0(\g)^\perp$. We then endow $GL_0(\g)^*$ with an Arens-type
product which contains $M(\g)$ as a closed subalgebra and
$M_a(\g)$ as a closed ideal which is a solid set with respect to
absolute continuity in $GL_0(\g)^*$. Among other things, we prove
that $\g$ is compact if and only if $GL_0(\g)^*$ has a non-zero
left (weakly) completely continuous element.
\\

\noindent{\bf Mathematics Subject Classification (2010).} 43A10, 43A15, 43A20, 47B07.

\noindent{\bf Key words.} Measure algebra, generalised functions
vanishing at infinity,
 introverted subspace, topological invariant mean,
completely continuous element.
\end{abstract}


\section{\large\bf Introduction}

All over this paper $\g$ is a locally compact group
with left Haar measure $\lambda$ and identity element $e$, and the notations
$C_c(\g)$ and $C_0(\g)$ refer to the space of all functions
with compact support and the space of
all functions vanishing at infinity, respectively. Moreover, the letter $M(\g)$ means that the measure algebra of $\g$ consisting of all complex regular Borel
measures on $\g$ with the total variation norm and the
convolution product ``$\ast$" defined by the formula
 $$
 \langle\mu\ast\nu ,g\rangle=\int_\g\int_\g g(xy)\;
 d\mu(x)\; d\nu(y)
 $$
 for all $\mu, \nu\in M(\g)$ and $g\in C_0(\g)$. It is folklore that $M(\g)$
is the first dual space of $C_0(\g)$ for the pairing
$$ \langle \mu,g \rangle:=\int_\g  g(x)\;d\mu(x)\quad\quad\quad\quad\quad
{\Big(}\mu\in M(\g),~g\in C_0(\g){\Big)}.$$


In the last thirty years, research on the second duals of Banach
algebras has mostly centred around the Banach algebras related to
locally compact groups, and has been dealt with by Lau {\it et
al.} in the works \cite{DLS2,DLS,GL,LP}. In particular, \cite{GL}
is the first important work devoted to the study of the second
duals of measure algebras. Among other things, the authors of
\cite{GL} have conjectured that the Banach algebra $M(\g)$ is
strongly Arens irregular, and its second dual, $M(\g)^{**}$,
determines $\g$ in the category of all locally compact groups.
Later on, the second duals of the measure algebras has been
studied in a series of papers. In particular, \cite{DLS} is the
second important work devoted to the study of $M(\g)^{**}$, where
most of the known results about these Banach algebras up to the
year 2012 can be found in it. Recall from \cite{DLS} that
$M(\g)^*$, the first dual space of $M(\g)$, as the second dual of
the C$^*$-algebra $C_0(\g)$ is a commutative unital C$^*$-algebra,
and therefore if $\widetilde{\g}$ denotes the hyper-Stonean
envelope of $\g$, then we can recognizing $M(\g)^*$ as
$C(\widetilde{\g})$, the space of all bounded complex-valued
continuous functions on $\widetilde{\g}$. It follows that
$M(\g)^{**}\cong M(\widetilde{\g})$, where $\cong$ denotes the
isometric algebra isomorphism and many authors, up to the year
2012, have used a type of this identification as a tool for the
study of $M(\g)^{**}$, see \cite{DLS} and the references therein
for more details.

Recently in the works \cite{EJN,JN2}, we studied the first and
second duals of measure algebras by the use of the theory of
generalised functions which have been introduced and investigated
by \u{S}re\u{i}dr \cite{S} and Wong \cite{W2,W3}. In those papers
we observed that $GL_0(\g)$, the space of all generalised
functions that vanishes at infinity, plays a crucial role in our
investigation. Motivated by this, here we study the C$^*$-algebra
$GL_0(\g)$ as an introverted subspace of $M(\g)^*$. In particular,
in the case where $\g$ is non-compact we show that any topological
left invariant mean on $GL(\g)$ lies in $GL_0(\g)^\perp$ which
demonstrates that the weak$^*$-closed subspace $GL_0(\g)^\perp$ of
$M(\g)^*$ is far from devoid of interest. We then endow
$GL_0(\g)^*$ with an Arens-type product which contains $M(\g)$ and
$M_a(\g)$ as a closed subalgebra and a closed ideal, respectively.
Among other things, we prove that the existence of a non-zero left
(weakly) completely continuous element in $GL_0(\g)^*$ is
equivalent to the compactness of $\g$.


\section{\large\bf Generalised functions: an overview}

In this section, we give a brief overview of generalised functions in the sense of
Wong \cite{W2}. Nevertheless, we shall require some facts
about the theory of C$^*$-algebra. For background on this theory,
we use \cite{MURP} as a reference and adopt that book's notation. Moreover, our notation and terminology are standard and, concerning Banach algebras related to locally compact groups, they
are in general those of the book \cite{HR} of Hewitt and Ross.
The reader will remark that this section is mostly taken from the papers by  Wong \cite{W2,W3}; this is because of, we actually need this section for the convenience of citation and a better exposition.

For any complex regular Borel measure $\mu$ on $\g$,
let $L^\infty(|\mu|)$ denote the Banach space of all essentially
bounded $|\mu|$--measurable complex functions $f_\mu$ on
${\mathcal G}$ with the essential supremum norm
$$\|f_\mu\|_{\mu,\infty}=\inf\Big{\{}\alpha\geq 0: |f_\mu|\leq\alpha,
~|\mu|-a.e.\Big{\}}.$$
Consider the product linear space
$\prod\{L^{\infty}(|\mu|):\mu\in M({\mathcal G})\}$. An element
$f=(f_\mu)_{\mu\in M({\mathcal G})}$ in this product is called a
{\it generalised function} if $f_\mu=f_\nu~|\mu|-a.e.$ for any
$\mu,\nu \in M({\mathcal G})$ with $\mu \ll \nu$,
where $\mu \ll \nu$ means that $|\mu|$
is absolutely continuous with respect to $|\nu|$. We note that
this condition implies that for given generalised functions $f=(f_\mu)_{\mu\in M(\g)}$  $$\sup
{\Big\{}\|f_\mu\|_{\mu,\infty}:\mu\in M(\g){\Big\}}<\infty;$$
otherwise, there
is a sequence $(\mu_n)$ in $M(\g)$ for which
$\|f_{\mu_n}\|_{{\mu_n},\infty}\geq n$ for all $n\in {\Bbb N}$.
Set $\mu=\sum_{n=1}^\infty 2^{-n}\|\mu_n\|^{-1}|\mu_n|$.
Then $\mu_n\ll\mu$, and hence
$\|f_\mu\|_{{\mu},\infty}\geq\|f_{\mu_n}\|_{{\mu_n},\infty}\geq n
$ for all $n\in{\Bbb N}$ which is a contradiction.


Now, following Wong \cite{W2} we use the notation $GL({\mathcal G})$ to denote the
commutative unital C$^*$-algebra of all generalised functions
endowed with the coordinatewise operations, the involution
$f\mapsto f^*$, where $f^*:=(\,\overline{f_\mu}\,)_{\mu\in M({\mathcal G})}$,
and norm $$\|f\|_\infty:=\sup{\Big\{}\|f_\mu\|_{\mu,\infty}:\mu\in
M(\g){\Big\}},$$ where $f=(f_\mu)_{\mu\in M(\g)}$ is in $GL(\g)$. The identity
element of $GL(\g)$ is of course the generalised function ${\bf 1}:=(1_\mu)_{\mu\in
M(\g)}$, where $1_\mu$ is the identity element of
$L^\infty(|\mu|)$. Moreover, we write $f=(f_\mu)_{\mu\in M(\g)}\geq 0$
to mean that the generalised function $f$ is positive in the C$^*$-algebra sense
and denote by $GL(\g)^+$ the set of all positive elements of $GL(\g)$.

\begin{remark}\label{positive}
{\rm It is not hard to check that a generalised function
$f=(f_\mu)_{\mu\in M(\g)}$ is positive in the C$^*$-algebra sense
if and only if $f_\mu(x)\geq 0$ for all $x\in\g$ and all $\mu\in
M(\g)$; see \cite[Page. 45]{MURP} and \cite[Page. 85]{W2} for more
information.}
\end{remark}


As a main result, Wong \cite{W2} has shown that for
each $f=(f_\mu)_{\mu\in M(\g)}$ in $GL(\g)$, the
equation
$$\big<\Psi(f),\zeta\big>:=\int_{\mathcal G}f_\zeta(x)\;d\zeta(x)\quad
\quad\left( \zeta\in M(\g) \right)$$ defines a linear functional
$\Psi(f)$ on $M(\g)$. In particular, the map $f\mapsto\Psi(f)$ is an isometric linear
mapping from $GL({\mathcal G})$ onto $M({\mathcal G})^*$; see
{\cite[Theorem 2.1 and Theorem 2.2]{W2}} and \cite{S} for the same
result in the special case where $\g$ is a certain locally compact abelian
group. In particular, any $L\in M({\mathcal G})^*$ can be
considered as a generalised function $\Psi^{-1}(L)$, and we do not
distinguish between a generalised function $f$ and its unique
corresponding linear functional $\Psi(f)$.
In particular, this duality allows us to consider $GL(\g)$ as a Banach
$M(\g)$-bimodule. In details, if $\zeta$ and $f=(f_\mu)_{\mu\in M(\g)}$ are arbitrary elements
of $M(\g)$ and $GL(\g)$, respectively. Then, one can considered the linear functionals $f\zeta$ and
$\zeta f$ on $M(\g)$ defined by
$$\big<f\zeta,\mu\big>=\big<f,\zeta\ast\mu\big>,\quad \big<\zeta f,\mu\big>=
\big<f,\mu\ast\zeta\big>\quad\quad\quad\quad\quad{\Big(}\mu\in M(\g){\Big)}.$$ In order to find the
generalised functions corresponding to these linear functionals,
following Wong \cite{W2} and \cite[Page 610]{W3}, we define
$\zeta\circ f\in\prod\{L^\infty (|\mu|):~\mu\in M(\g)\}$ as
$$(\zeta\circ f)_\mu=l_\zeta f_{\zeta\ast\mu}\quad\quad(\mu\in
M(\g)),$$
where
\begin{eqnarray*}
l_\zeta f_{\zeta\ast\mu}(y)=\int_\g
f_{\zeta\ast\mu}(xy)\;d\zeta(x)\quad\quad {\hbox{for}}
~~~|\mu|-a.e ~~y\in\g.
\end{eqnarray*}
Then $\zeta\circ f$ is again a
generalised function such that
\begin{eqnarray*}
\big<\zeta\circ f,\mu \big>=\big<f,\zeta\ast\mu\big>\quad\quad (\mu\in M(\g));
\end{eqnarray*}
see pages 88 and 89 of \cite{W2}.
So, $\zeta\circ f$ is the generalised function corresponding to
the functional $f\zeta\in M(\g)^*$ such that $\Psi(\zeta\circ f)=f\zeta$.
Also, by using the right convolution notation, we can show that $f\circ\zeta$ is the
generalised function corresponding to the functional $\zeta f\in
M(\g)^*$. In what follows,
we do not distinguish between the linear functionals $f\zeta$ and $\zeta f$ and their
corresponding generalised functions.


Later we will need the following remark in our present investigation.

\begin{remark}\label{tojih}
{\rm Suppose that $BM(\g)$ denote the Banach space of all
bounded Borel measurable functions on $\g$ with the supremum norm
$\|\cdot\|_u$. Then each $f\in BM(\g)$ may be regarded as an element
$(f_\mu)_{\mu\in M(\g)}$ in $GL(\g)$, where for each $\mu\in
M(\g)$ the functions $f_\mu$ denotes the
equivalent class of $f$ in $L^\infty(|\mu|)$. Hence $BM(\g)$ can be considered as a closed subspace of
$GL(\g)$ containing the space $C_b(\g)$ of all complex-valued continuous bounded
functions on $\g$. Moreover, each $f\in BM(\g)$ may be regarded as an
element in $M(\g)^*$ by the pairing
$\big<f,\mu\big>=\int_\g f\;d\mu$ ($\mu\in M(\g)$). In this case, the restriction of the map $\Psi$ to $BM(\g)$
is precisely the embedding of $BM(\g)$ into $M(\g)^*$.
}
\end{remark}


\section{\large\bf Generalised functions that vanish at infinity}


We commence this section by recalling the main object of the work
which is introduced and studied by the authors in \cite{JN2}.

\begin{definition}
A generalised function $f=(f_\mu)_{\mu\in M(\g)}$ {\it vanishes at
infinity} if for each $\varepsilon>0$, there is a compact subset
$K_\varepsilon$ of ${\mathcal G}~$ for which
$\|f_\mu\chi_{{\mathcal G}\setminus
K_\varepsilon}\|_{\mu,\infty}<\varepsilon$ for all $\mu\in
M({\mathcal G})$; formally
\begin{eqnarray*}
\forall\varepsilon>0\quad\exists K_\varepsilon\in{\mathcal
K}(\g)\quad{\hbox{s.t.}}\quad \forall\mu\in
M(\g),\quad|f_\mu(x)|<\varepsilon\quad{\hbox{for}}~|\mu|-{\hbox{almost~all}}~x\in
\g\setminus K_\varepsilon,
\end{eqnarray*}
where $\chi_{K_\varepsilon}$ denotes the characteristic function
of $K_\varepsilon$ on ${\mathcal G}$ and ${\mathcal K}(\g)$
denotes the set of all compact subsets in $\g$.
\end{definition}

We denote by $GL_0({\mathcal G})$ the C$^*$-subalgebra of $GL(\g)$
consisting of all generalised functions that vanish at infinity.

The aim of the present section is to study some aspects of
$GL_0(\g)$ as a C$^*$-subalgebra of $GL(\g)$. Let us give a simple
but important result whose proof involves nothing more than
routine calculations.

\begin{lemma}\label{appiden}
Suppose that ${\mathcal K}(\g)$
is directed downwards and for each $\alpha$,
$u_{K_\alpha}\in C_c(\g)$ is
chosen such that $0\leq u_{K_\alpha}\leq 1$ and
$u_{K_\alpha}(x)=1$ for all $x\in K_\alpha$. Then $(u_{K_\alpha})$
is a bounded approximate identity for
$GL_0(\g)$.
\end{lemma}


Our next result shows that the subspaces
$$M(\g)\circ GL_0(\g):={\Big\{}\zeta\circ f:~\zeta\in M(\g)~{\hbox{and}}
~f\in GL_0(\g){\Big\}}$$
and
$$ GL_0(\g)\circ M(\g):={\Big\{}f\circ\zeta:~\zeta\in M(\g)
~{\hbox{and}}~f\in GL_0(\g){\Big\}}$$
of $GL(\g)$ coincide with $GL_0(\g)$.

\begin{lemma}\label{j4}
The following assertions hold.
\newcounter{j1215}
\begin{list}%
{\rm(\roman{j1215})}{\usecounter{j1215}}
\item $M(\g)\circ GL_0(\g)=GL_0(\g)$,
\item $GL_0(\g)\circ M(\g)=GL_0(\g)$.
\end{list}
\end{lemma}
{\noindent Proof}. We prove the first; the proof of the second is
similar.  Since the inclusion $GL_0(\g)\subseteq \delta_e\circ
GL_0(\g)\subseteq M(\g)\circ GL_0(\g)$, it will be enough to prove
the reverse inclusion. To this end, let $\zeta\in M(\g)$,
$f=(f_\mu)_{\mu\in M(\g)}\in GL_0(\g)$ and $\epsilon>0$ be given.
Without loss of generality, we may assume that $\zeta$ is non-zero
and positive and that $f\neq 0$. By the regularity of $\zeta$, we
can choose a compact subset $K_1$ of $\g$ such that
$0<\zeta(\g\setminus K_1)<\epsilon/2\|f\|_\infty$. Also, since $f$
vanishes at infinity, there is a compact subset $K_2$ in $\g$ with
$\|f-\chi_{K_2}f\|_\infty<\epsilon/2\|\zeta\|$. Therefore
\begin{equation}\label{E1}
\|\zeta\circ f-(\chi_{K_1}\zeta)\circ(\chi_{K_2}f)\|_\infty\leq
\|\zeta-\chi_{K_1}\zeta\|\;\|f\|_\infty+\|\chi_{K_1}\zeta\|\|f-\chi_{K_2}f\|_\infty
\leq \epsilon,
\end{equation}
where $\chi_{K_1}\zeta$ is the measure in $M(\g)$ defined on each Borel
subsets $A$ of $\g$ by
$\chi_{K_1}\zeta(A)=\int_A\chi_{K_1}\;d\zeta$. Now suppose that $\mu$ is an arbitrary element
of $M(\g)$. Observe that
\begin{eqnarray*}
\Big{(} (\chi_{K_1}\zeta)\circ(\chi_{K_2}f)
\Big{)}_\mu(x)=\int_{K_1}
\chi_{K_2}(yx)f_{(\chi_{K_1}\zeta)\ast\mu}(yx)\;d\zeta(y).
\end{eqnarray*}
For each $x\in\g\setminus
{K_1}^{-1}K_2$, we get $K_1x\subseteq\g\setminus
K_2$, and hence
$\Big{(} (\chi_{K_1}\zeta)\circ(\chi_{K_2}f) \Big{)}(x)=0$
for $\mu$--almost all $x\in\g\setminus {K_1}^{-1}K_2$. Thus, inequality (\ref{E1}) implies that
$$|(\zeta\circ f)_{\mu}(x)|<\varepsilon,\quad\quad{\hbox{
for}}~ \mu-{\hbox{almost~ all}}~ x\in\g\setminus {K_1}^{-1}K_2.$$ It follows that $\zeta\circ f\in GL_0(\g)$.
We have now completed the proof of the lemma.
$\hfill\blacksquare$\\


Now, let $M_a(\g)$ be the closed ideal of $M(\g)$ consisting of
all absolutely continuous measures with respect to $\lambda$ and
let $L^1(\g)$ denote the group algebra of $\g$ as defined in
{\cite[Theorem 14.17 and 14.18]{HR}}. Then, the Radon-Nikodym
Theorem can be interpreted as an identification  of $M_a(\g)$ with
$\{\nu_\varphi:~\varphi\in L^1(\g)\}$, where $\nu_\varphi$ is the
measure in $M(\g)$ defined on each Borel subset $A$ of $\g$ by
$\nu_\varphi(A)=\int_A\varphi\;d\lambda$. This allows us to show
that $M_a(\g)^*$, the first dual space of $M_a(\g)$, is
$L^\infty(\g)$, where $L^\infty(\g)$ denotes the Lebesgue space as
defined in {\cite[Definition 12.11]{HR}} equipped with the
essential supremum norm. Given any $\sigma\in M_a(\g)$ and $g\in
L^\infty(\g)$, define the complex-valued functions $g\star\sigma$
and $\sigma\star g$ on $\g$ by
$$(g\star\sigma)(x)=\big<g,\delta_x\ast\sigma\big>=\int_\g g(xy)\;d\sigma(y)$$
and
$$(\sigma\star g)(x)=\big<g,\sigma\ast\delta_x\big>=\int_\g g(yx)\;d\sigma(y)$$
for all $x\in\g$, where $\delta_x$ denotes the Dirac measure at
$x$. Then, it is not hard to check that the functions
$g\star\sigma$ and $\sigma\star g$ are in $C_b(\g)$; this is
because of, $M_a(\g)$ can be identified with all $\nu\in M(\g)$
such that the map $x\mapsto\delta_x\ast|\nu|$ and $x\mapsto
|\nu|\ast\delta_x$ from $\g$ into $M(\g)$ are norm continuous, see
for example \cite[19.27 and 20.31 ]{HR}. In particular,  if
${\mathcal{P}}:GL(\g)\rightarrow L^\infty(\g)$ is the adjoint of
the natural embedding from $M_a(\g)$ into $M(\g)$, then
$\mathcal{P}$ is the restriction mapping and hence norm decreasing
and onto.

For the formulation of the following statements we recall Remark
\ref{tojih} which allows us to consider $C_b(\g)$ as a closed
subspace of $GL(\g)$ containing $C_0(\g)$.

\begin{lemma}\label{j11}
If $\sigma$ is an arbitrary element of $M_a(\g)$, then for given
$f=(f_\mu)_{\mu\in M(\g)}$ in $GL(\g)$ and all $\mu\in M(\g)$, we
have
\newcounter{j1216}
\begin{list}%
{\rm(\roman{j1216})}{\usecounter{j1216}}
\item $(\sigma\circ f)_\mu=\sigma\star{\mathcal P}(f)$,
$|\mu|$--a.e.,
\item $(f\circ\sigma)_\mu={\mathcal P}(f)\star\sigma$,
$|\mu|$--a.e.
\end{list}
In particular, $(\sigma\circ f)_\mu$ and
$(f\circ\sigma)_\mu$ are in $C_b(\g)$ for all $\mu\in M(\g)$.
\end{lemma}
{\noindent Proof}. We prove the assertion (i); the proof of (ii)
is similar. First note that $\sigma\circ f$ is the generalised
function $h=(h_\mu)_{\mu\in M(\g)}$, where for each $\mu\in M(\g)$
satisfies the following equality
\begin{eqnarray*}
h_\mu(x)=(\sigma\circ f)_\mu(x)=l_\sigma f_{\sigma\ast\mu}(x)
=\int_\g f_{\sigma\ast\mu}(yx)\;d\sigma(y).
\end{eqnarray*}
On the
other hand, for an arbitrary $\mu$ in $M(\g)$ and any Borel subset
$A$ of $\g$, since $\chi_A\mu\ll\mu$, we have
\begin{eqnarray*}
\int_\g \chi_Ah_\mu\;d\mu&=&
\int_\g f_{\sigma\ast\chi_A\mu}\;d(\sigma\ast\chi_A\mu)\\
&=&\big<{\mathcal P}(f),\sigma\ast\chi_A\mu\big>\\
&=&\int_\g\int_\g
{\mathcal P}(f)(yx)\;d\sigma(y)d(\chi_A\mu)(x)\\
&=&\int_\g\chi_A\;(\sigma\star{\mathcal P}(f))\;d\mu.
\end{eqnarray*}
 Hence, $(\sigma\circ f)_\mu=\sigma\star{\mathcal
P}(f)$ $|\mu|$--a.e. ($\mu\in M(\g)$). It follows that
$\sigma\circ f\in C_b(\g)$.$\hfill\blacksquare$\\


Now, in light of Lemmas \ref{j4} and \ref{j11}, the following
proposition is now immediate.

\begin{proposition}\label{05}
The following assertions hold.
\newcounter{j1217}
\begin{list}%
{\rm(\roman{j1217})}{\usecounter{j1217}}
\item $M_a(\g)\circ GL_0(\g)=C_0(\g)$,
\item $GL_0(\g)\circ M_a(\g)=C_0(\g)$.
\end{list}
\end{proposition}


Recall from \cite[Page. 90]{W2} that a linear functional $\m$ in
$GL(\g)^*$ is called a { mean} if $\m({\bf 1}) = 1$ and $\m(f)\geq
0$ whenever $f\in GL(\g)$ with $f\geq 0$, and it is { topological
left invariant} if $\m(\zeta\circ f)=\m(f)$ for all $f\in GL(\g)$
and $$\zeta\in P(\g)={\Big\{}\nu\in M(\g):~\nu\geq
0~{\hbox{and}}~\|\nu\|=1{\Big\}}.$$ In \cite[Theorem 4.1]{W2},
Wong proved that $GL(\g)$ has a topological left invariant mean if
and only if $M(\g)^*$ has a topological left invariant mean. In
particular, he showed that $\Psi^*$, the adjoint of $\Psi$, maps
the set of all topological left invariant means on $M(\g)^{*}$
onto that of $GL(\g)$. Related to this result, we have the
following result which asserts that in the case where $\g$ is
non-compact, then any topological left invariant mean on $GL(\g)$
lies in $GL_0(\g)^\perp$, where here and in the sequel,
$GL_0(\g)^\perp$ denotes the following weak$^*$-closed subspace of
$GL(\g)^*$
$${\Big\{}\m\in GL(\g)^*:~\big<\m,f\big>=0~
{\hbox{for~all~}}f\in GL_0(\g){\Big\}}.$$ In fact, the next
result shows that $GL_0(\g)^\perp$ is far from devoid of interest.

\begin{proposition}\label{mean}
If $\g$ is non-compact, then
any topological left invariant mean on $GL(\g)$ lies in $GL_0(\g)^\perp$.
\end{proposition}
{\noindent Proof.} Suppose that $\m$ is a topological left
invariant mean on $GL(\g)$. First note that the non-compactness of $\g$ implies
that there exists a sequence $(x_n)$ of
disjoint elements of $\g$ and a compact symmetric neighborhood $V$
of $e$ such that the sets $x_nV$ for all ${n\in{\Bbb N}}$ are pairwise
disjoint; see {\cite[11.43(e)]{HR}}. Now, it is not hard to check
that $\chi_{xV}=\delta_{{x}^{-1}}\circ\chi_V$, $|\mu|-$a.e. for
all $x\in\g$ and $\mu\in M(\g)$. Moreover, by Remark
\ref{tojih}, for each $p\in{\Bbb N}$ the function
$\sum_{n=1}^p\chi_{x_nV}=\sum_{n=1}^p\delta_{{x_n}^{-1}}\circ\chi_V$
is in $GL(\g)$ for which
$\sum_{n=1}^p\chi_{x_nV}\leq {\bf 1}$. It follows that
$$
p\big<\m,\chi_{V}\big>=\big<\m,\sum_{n=1}^p\chi_{x_nV}\big>\leq 1\quad\quad\quad(p\in{\Bbb N}).
$$
Thus $\big<\m,\chi_V\big>=0$ and thus we have
\begin{eqnarray}\label{41194}
\big<\m,\chi_{xV}\big>=\big<\m,\delta_{{x}^{-1}}\circ\chi_{V}\big>=0\quad\quad\quad(x\in\g).
\end{eqnarray}

Now suppose that $f=(f_\mu)_{\mu\in M(\g)}$ is a non-zero element
of $GL_0(\g)$. The proof will be completed by showing that
$\big<\m,f\big>=0$. To this end, without loss of generality, we
may assume that $\|f\|_\infty=1$. Then, since
$f^*f=(|f_\mu|^2)_{\mu\in M(\g)}$ vanishes at infinity, for given
$\varepsilon>0$ one can choose $y_1,...,y_q\in\g$ such that
$$
|f_\mu|^2\leq\sum_{i=1}^q\chi_{y_iV}+
\varepsilon\quad|\mu|-{\hbox{a.e.}}\quad\quad\quad\quad{\Big(}\mu\in
M(\g){\Big)}.
$$
Now, by considering $h=\sum_{i=1}^q\chi_{y_iV}+\varepsilon\in
BM(\g)$ as an element of $GL(\g)$, we have $f^*f\leq h$, see
Remark \ref{positive}. Hence, in light of \cite[Theorem
3.3.2]{MURP} and equality (\ref{41194}), we see that
$$|\big<\m,f\big>|^2\leq\big<\m,f^*f\big>\leq\big<\m,h\big>\leq\varepsilon.$$
It follows that
$\big<\m,f\big>=0$. Hence,
$\m\in GL_0(\g)^\perp$.$\hfill\blacksquare$\\


\section{\large\bf $GL_0(\g)^*$  as a subalgebra of  $M(\g)^{**}$}

As we know, there exists two natural products on $M(\g)^{**}$ extending the
one on $M(\g)$, known as the first and second Arens products of
$M(\g)^{**}$. The first Arens product on $M(\g)^{**}$ is defined
in three steps as follows. For $\m$, $\n$ in $M(\g)^{**}$, the
element $\m\odot \n$ of $M(\g)^{**}$ is defined by
$$\big<\m\odot \n,f\big>=\big<\m,\n f\big>\quad\quad (f\in GL(\g)),$$
where $\big<\m f,\zeta\big>=\big<\m,f\zeta\big>$ and
$f\zeta=\zeta\circ f$ for all $\zeta\in M(\g)$. Equipped with this
product, $M(\g)^{**}$ is a Banach algebra which contains $M(\g)$
as a subalgebra. Moreover, by the duality relation between $M(\g)$
and $GL(\g)$, there exists a unique generalised function $h\in
GL(\g)$ such that $\m f=\Psi(h)$; In what follows, we denote the
generalised function $\Psi^{-1}(\m f)$ corresponding to $\m f\in
M(\g)^*$, by $\m f$. Moreover, $M(\g)$ and $GL_0(\g)$ are in
duality with respect to the natural bilinear map given for each
$\eta\in M(\g)$ and $f=(f_\mu)_{\mu\in M(\g)}$ in $GL_0(\g)$ by
$\big<f,\eta\big>=\int_\g f_\eta\;d\eta.$ Therefore, $M(\g)$ may
be identified with a closed subspace of $GL_0(\g)^*$. Furthermore,
if $f=(f_\mu)_{\mu\in M(\g)}\in GL_0(\g)$ and $\zeta\in M(\g)$,
then, by Lemma \ref{j4}, $\zeta\circ f$ and $f\circ\zeta$ are also
in $GL_0(\g)$ and
$$\big<\zeta\circ f,\nu\big>=\big<f,\zeta\ast\nu\big>, \quad\quad{\hbox{and}}~
\quad\big<f\circ\zeta,\nu\big>=\big<f,\nu\ast\zeta\big>,$$
for all $\nu\in M(\g)$. Hence the product ``$\odot$" is well defined on $GL_0(\g)^*$ and $GL_0(\g)^*$
is a Banach algebra with this product,
if we show that $GL_0(\g)$ is a topologically introverted subspace of
$GL(\g)$. To this end, we have the following result.


\begin{proposition}\label{st}
The space $GL_0(\g)$ is left {\rm(}right{\rm)} topologically
introverted in $GL(\g)$; that is, ${\rm\m} f\in GL_0(\g)$
{\rm(}$f{\rm\m}\in GL_0(\g)${\rm)} for all ${\rm\m}\in GL_0(\g)^*$
and $f\in GL_0(\g)$.
\end{proposition}
{\noindent Proof}. We only show that $GL_0(\g)$ is a left
topologically introverted subspace of $GL(\g)$; the proof of the other assertion is similar.
 To this end, let $\m\in GL_0(\g)^*$, $f=(f_\mu)_{\mu\in
M(\g)}\in GL_0(\g)$ and $\varepsilon>0$ be given. Since
$GL_0(\g)$ is spanned by its positive elements, we can suppose that
$\m\geq 0$.
Also, since $f$ vanish at infinity, there is a compact set $B$ in
$\g$ with $|f_\mu(x)|<\varepsilon$ for $\mu$-almost all
$x\in\g\setminus B$ $(\mu\in M(\g))$.

Now let $\varrho$ denote the restriction
of $\m$ to $C_0(\g)$. Then there exists a compact subset $K$
of $\g$ such that $\varrho(\g\setminus K)<\varepsilon/2$.
In particular, if $\m_K$ denote the continuous linear functional on
$GL_0(\g)$ defined by
$$\big<\m_K,h\big>:=\big<\m,h-u_Kh\big>\quad\quad\quad{\Big(}h\in GL_0(\g){\Big)},$$
where $u_{K}$ is a fixed function in $C_c(\g)$ such that $0\leq u_{K}\leq 1$, and
$u_{K}(x)=1$ for all $x\in K$. Then the positivity of the linear
functional $\m_K$ on $GL_0(\g)$ implies that
$\|\m_K\|=\lim_\alpha\big<\m_K,u_{K_\alpha}\big>$, where $(u_{K_\alpha})$ is the net introduced in Lemma \ref{appiden}. Hence, there exists
$\alpha_0$ such that
$$\|\m_K\|-\frac{\varepsilon}{2}\leq\big<\m_K,u_{{K_{\alpha_0}}}\big>\leq\|\m_K|_{C_0(\g)}\|.$$
It follows that $\|\m_K\|<\varepsilon$. Indeed,
\begin{eqnarray*}
\|\m_K|_{C_0(\g)}\|
&=&\sup{\Big\{}|\big<\varrho,g-u_Kg\big>|:~g\in C_0(\g)~{\rm{and}}~\|g\|\leq 1{\Big\}}\\
&=&\sup{\Big\{}{\Big|}\big<\chi_{\g\setminus K}\varrho,(g-u_Kg)\big>{\Big|}:~g\in C_0(\g)~{\rm{and}}~\|g\|\leq 1{\Big\}}\\
&\leq&\|\chi_{\g\setminus K}\varrho\|\\&=&\varrho(\g\setminus K).
\end{eqnarray*}

If now, $\nu$ is an arbitrary probability measure in $M(\g)$, then
$\zeta:=(\chi_{\g\setminus BK^{-1}})\nu$ is a measure in $M(\g)$ for which ${\hbox
{supp}}(\zeta)\subseteq\g\setminus BK^{-1}$. Further, choose a compact
subset $D$ in $\g$ for which $D\subseteq\g\setminus BK^{-1}$ and that
$|\zeta|(\g\setminus D)<\varepsilon$. Trivially, for each
$x\in\g\setminus D^{-1}B$, we see that $Dx\subseteq\g\setminus B$, and
therefore for each $\mu\in M(\g)$ we have
\begin{eqnarray*}
|(\zeta\circ f)_\mu(x)|\leq\int_{\g\setminus D} |f_{\zeta\ast\mu}(yx)|\;d|\zeta|(y)
+\int_{D} |f_{\zeta\ast\mu}(yx)|\;d|\zeta|(y)
\leq\varepsilon(\|f\|_\infty+1);
\end{eqnarray*}
that is,
$|(\zeta\circ f)_\mu(x)|\leq\varepsilon(\|f\|_\infty+1)$
for $\mu$-almost all $x\in \g\setminus D^{-1}B$.
In particular, since $D^{-1}B\cap K=\emptyset$, we see that
$$\|(\zeta\circ f)\chi_K\|_\infty=\sup_{\mu\in M(\g)}
\|(\zeta\circ f)_\mu\chi_K\|_{\mu,\infty}\leq\varepsilon
(\|f\|_\infty+1).$$
Thus
\begin{eqnarray*}
\int_{\g\setminus BK^{-1}}(\m f)_\zeta(x)\;d\zeta(x)&=&\big<\m f,\zeta\big>\\
&=&\big<\m,u_K(\zeta\circ f)\big>+\big<\m_K,\zeta\circ f\big>\\
&\leq&\varepsilon(\|f\|_\infty+1)\|\m\|+\varepsilon\|\zeta\|\|f\|_\infty.
\end{eqnarray*}
On the other hand, since $\zeta\ll\nu$, we have
\begin{eqnarray*}
\int_{\g\setminus
BK^{-1}}(\m f)_\zeta(x)\;d\zeta(x)=\int_{\g\setminus
BK^{-1}}(\m f)_\nu(x)\;d\zeta(x)
=\int_{\g\setminus BK^{-1}}(\m f)_\nu(x)\;d\nu(x).
\end{eqnarray*}
This shows that, if $\nu\in M(\g)$, then
$(\m f)_\nu(x)\leq\varepsilon[(\|f\|_\infty+1)\|\m\|+\|f\|_\infty]$ for
$\nu$-almost all $x\in\g\setminus BK^{-1}$, and thus $\m f\in
GL_0(\g)$. $\hfill\blacksquare$\\


A linear functional $\m$ in
$GL_0(\g)^{*}$ (resp. $M(\g)^{**}$) has \emph{compact
carrier} if there exists a compact set $K$ in $\g$ such that $
\langle \m,f \rangle=\langle \m , \chi_{K}f \rangle$ for all $f\in
GL_{0}(\g)$ (resp. $GL(\g)$);
such a compact set $K$ is called a compact carrier for $\m$.
In the sequel, the notation $M_c(\g)^{**}$ is used to denote the norm
closure of functionals in $M(\g)^{**}$ with
compact carrier.


Now, with an argument similar to the proof of Propositions 2.6, 2.7 and Theorems 2.8 and 2.11 in \cite{LP}, one can prove the following result which in particular shows that the restriction map is an isometric algebra isomorphism from $M_c(\g)^{**}$ onto
$GL_0({\mathcal G})^*$. In other word, this result allows us to view $GL_0(\g)^*$
as a subalgebra of $M(\g)^{**}$.

\begin{theorem}\label{02}
The following
assertions hold.
\newcounter{j1212}
\begin{list}%
{\rm(\roman{j1212})}{\usecounter{j1212}}
\item Functionals in $GL_0(\g)^*$ with compact
carriers are norm dense in $GL_0(\g)^*$.
\item If $\m$ and $\n$ are elements
in $GL_0(\g)^*$ {\rm(}resp. $M_c(\g)^{**}${\rm)} with
compact carriers $K$ and $K'$ respectively, then $\m\odot \n$
has compact carrier $KK'$.
\item The restriction map is an isometry and an algebra isomorphism from
$M_c(\g)^{**}$ onto $GL_0(\g)^{*}$.
\item $M(\g)^{**}=GL_{0}(\g)^*\oplus
GL_{0}(\g)^\perp$.
In fact, any ${\rm\m}\in M(\g)^{**}$
has a unique decomposition
${\rm\m}={\rm\m}_*+{\rm\m}_\perp$,
where ${\rm\m}_*\in GL_{0}(\g)^*$,
${\rm\m}_\perp\in GL_{0}(\g)^\perp$ and
$\|{\rm\m}\|=\|{\rm\m}_*\|+\|{\rm\m}_\perp\|$. Moreover, ${\rm\m}\geq 0$ if and only if
${\rm\m}_*\geq 0$ and ${\rm\m}_\perp\geq 0$.
\item $GL_0(\g)^\perp$ is a
weak$^*$-closed ideal of $M(\g)^{**}$.
\item $GL_0(\g)^*$ is a left or right
ideal of $M(\g)^{**}$ if and only if $\g$ is compact.
\item $M(\g)$ is a left or right ideal
of $GL_0(\g)^*$ if and only if $\g$ is discrete.
\item $M_a(\g)$ is a two-sided ideal in $GL_0(\g)^*$.
\end{list}
\end{theorem}
{\noindent Proof}. The details are omitted and we only give the
proof for (ii) and (viii).

(ii) Suppose that $\m$ and $\n$ are elements
in $GL_0(\g)^*$ {\rm(}resp. $M_c(\g)^{**}${\rm)} with
compact carriers $K$ and $K'$, respectively, and that $f$ is an arbitrary element in
$GL_0(\g)$ (resp. $GL(\g)$). First, observe that
\begin{eqnarray*}
\big<\m\odot \n,f\big>=\big<\m,\n f\big>=\big<\m,\chi_{K}(\n f)\big>.
\end{eqnarray*}
On the one hand, for $\mu\in M(\g)$ we have
$\zeta:=\chi_{K}\mu\ll\mu$, and so
\begin{eqnarray*}
\big<\chi_{K}(\n f),\mu\big>=\int_{\g} \chi_{K}(\n f)_\mu \;d\mu
=\int_{\g} (\n f)_{\zeta}\;d\zeta
=\big<\n,\chi_{K'}(\zeta\circ f) \big>.
\end{eqnarray*}
On the other hand, $\chi_{K'}(\zeta\circ
f)=\zeta\circ(\chi_{KK'}f)$. Indeed, for each $\nu\in M(\g)$ and $h\in C_0(\g)$, we have
\begin{eqnarray*}
{\Big(}\chi_{KK'}(\zeta\ast\nu){\Big)}(h)&=&\int_{\g} \chi_{KK'}h \;d(\zeta\ast\nu)\\
&=&\int_{\g}\int_{\g} \chi_{KK'}(xy)\chi_{K}(x)h(xy)\;d\mu(x)\;d\nu(y)\\
&=&\int_{\g}\int_{\g} \chi_{K'}(y)\chi_{K}(x)h(xy)\;d\mu(x)\;d\nu(y)\\
&=&{\Big(}\zeta\ast(\chi_{K'}\nu){\Big)}(h),
\end{eqnarray*}
and this implies that
\begin{eqnarray*}
\big<\chi_{K'}(\zeta\circ f),\nu\big>=\big<\zeta\circ f,\chi_{K'}\nu\big>
=\big<f,\chi_{KK'}(\zeta\ast\nu)\big>
=\big<\zeta\circ(\chi_{KK'}f),\nu\big>.
\end{eqnarray*}
Hence by using these equalities, we have
\begin{eqnarray*}
\big<\chi_{K}(\n f),\mu\big>=\big<\n,\zeta\circ(\chi_{KK'}f)\big>
=\big<\n(\chi_{KK'}f),\zeta\big>
=\big<\chi_{K}(\n(\chi_{KK'}f)),\mu\big>.
\end{eqnarray*}
Consequently
\begin{eqnarray*}
\big<\m\odot
\n,f\big>=\big<\m,\chi_{K}(\n f)\big>
=\big<\m,\chi_{K}(\n(\chi_{KK'}f))\big>
=\big<\m,\n(\chi_{KK'}f)\big>
=\big<\m\odot \n,\chi_{KK'}f\big>.
\end{eqnarray*}
It follows that $\m\odot \n$ has compact carrier $KK'$.

(viii) That $M_a(\g)$ is a closed subalgebra of $GL_0(\g)^*$ is trivial.
Now, suppose that $\sigma\in M_a(\g)$ and $\m\in GL_0(\g)^*$. We
show that $\m\odot \sigma\in GL_0(\g)^*$; that $\sigma\odot \m\in
GL_0(\g)^*$ is similar. Let $\zeta$ denote the restriction of $\m$
to $C_0(\g)$. Since $M_a(\g)$ is an ideal in $M(\g)$, we have
$\zeta\ast\sigma\in M_a(\g)$.
We now invoke Proposition \ref{05} to conclude that
\begin{eqnarray*}
\big<\m\odot\sigma,f\big>=\big<\zeta,f\circ\sigma\big>=\big<\zeta\ast\sigma,f\big>
\quad\quad\quad{\Big(}f\in GL_0(\g){\Big)},
\end{eqnarray*}
whence $\m\odot \sigma=\zeta\ast\sigma\in M_a(\g)$.
$\hfill\blacksquare$\\


As usual, for a locally compact space $X$, we say that a subset
$S\subseteq M(X)$, the Banach space of all complex regular Borel
measures on $X$, is solid with respect to absolute continuity, if
$t\in S$ wherever $t\ll s$, for some $s\in S$. Now, as an
application of Theorem \ref{02} above, by a method similar to that
of \cite[Lemma 5 and Theorem 6]{GM}, one can obtain the following
generalization of that theorem; The reader will remark that the
compactness of $\g$ is assumed in that proof only to conclude that
$M_a(\g)$ is an ideal in $M(\g)^{**}$ whereas $M_a(\g)$ always is
an ideal of $GL_0(\g)^*$. The details are omitted.

\begin{theorem}\label{14}
$M_a(\g)$ is the unique
minimal proper closed subset of $GL_0(\g)^*$ which is an algebraic
ideal and a solid set with respect to absolute continuity in
$GL_0(\g)^*$.
\end{theorem}


Next we turn our attention to the study of left (weakly) completely
continuous elements of $GL_0(\g)^*$. To this end, recall that if $\cal A$ is a Banach
algebra, then $a\in{\cal A}$ is said to be a left (weakly) completely continuous element of ${\cal A}$ whenever the operator $\ell_a:b\mapsto ab$ is (weakly) compact operator on ${\cal A}$.

In what follows,
for $I\subseteq GL_0(\g)^*$, the left annihilator of $I$ is
denoted by ${\rm lan}(I)$ and defined by
$${\rm lan}(I)={\Big\{}\textsf{l}\in GL_0(\g)^*: \textsf{l}\odot I=\{0\}{\Big\}},$$
also, the right annihilator of $I$ is denoted by ${\rm ran}(I)$
and define by
$${\rm ran}(I)={\Big\{}\textsf{r}\in GL_0(\g)^*: I\odot \textsf{r}=\{0\}{\Big\}}.$$
Moreover, the letter
$L_0^\infty(\g)$ means that the C$^*$-subalgebra of $L^\infty(\g)$ consisting
of all functions $g$ on $\g$ such that for each $\varepsilon>0$,
there is a compact subset $K$ of $\g$ for which $|g(x)|<\varepsilon$ for all
$x\in\g\setminus K$.


\begin{theorem}
The following
assertions hold.
\newcounter{j1218}
\begin{list}%
{\rm(\roman{j1218})}{\usecounter{j1218}}
\item If $\sigma\in M_a(\g)$, then $\sigma$ is a left
{\rm(}weakly{\rm)} completely continuous element of $M_a(\g)$ if and only if
$\sigma$ is a left {\rm(}weakly{\rm)} completely continuous element of $GL_0(\g)^*$.
\item Any left {\rm(}weakly{\rm)} completely continuous
element $\rm\textsf{m}$ of $GL_0(\g)^*
$ has the form ${\rm\textsf{m}}=\sigma+{\rm\textsf{r}}$ for some $\sigma\in M_a(\g)$ and
${\rm\textsf{r}}\in{\rm ran}({\mathcal P}^*(L_0^\infty(\g)^*))$.
\end{list}
\end{theorem}
{\noindent Proof.} We only give the proof for
left completely continuous element.

(i) The direct implication being trivial, we give the proof of the backward implication
only. To this end, suppose that $\sigma\in M_a(\g)$ is a left
completely continuous element of $M_a(\g)$. Then,
the closure of the following set is compact in $M_a(\g)$
$${\Big\{}\sigma\ast\upsilon:
~\upsilon\in M_a(\g),~\|\upsilon\|\leq 1{\Big\}}.$$
On the other hand, if $(e_\alpha)_{\alpha}$ is an
approximate identity for $M_a(\g)$ bounded by one, then for each $\alpha$ and $\m\in GL_0(\g)^*$ with $\|\m\|\leq 1$,
we have
\begin{eqnarray*}
\|\sigma\odot \m-\sigma\ast(e_\alpha\odot \m)\|\leq\|\sigma-\sigma\ast e_\alpha\|.
\end{eqnarray*}
This together with the fact that $M_a(\g)$ is an ideal in $GL_0(\g)^*$ implies that
$$
{\Big\{}\sigma\odot \m:~\m\in GL_0(\g)^*,~\|\m\|\leq 1{\Big\}}\subseteq
\overline{{\Big\{}\sigma\ast\upsilon:
~\upsilon\in M_a(\g),~\|\upsilon\|\leq 1{\Big\}}}^{_{M_a(\g)}}.
$$
Thus the operator
$\ell_\sigma:GL_0(\g)^*\rightarrow GL_0(\g)^*$ is compact.

(ii) Suppose that $\m$ is a left completely continuous element of $GL_0(\g)^*$. Then,
since $M_a(\g)$ is an ideal in $GL_0(\g)^*$, the operator
$\ell_\m|_{M_a(\g)}$ is a compact operator on $M_a(\g)$. From this,
we can conclude that there exists
$\sigma\in M_a(\g)$ such that $\ell_\m=\ell_\sigma$ on $M_a(\g)$, see \cite{ake}.
In particular, Proposition \ref{05} implies that
$\big<\m,f\big>=\big<\sigma,f\big>$
for all $f\in C_0(\g)$, and thus we have
$\upsilon\odot\m=\upsilon\ast\sigma$ ($\upsilon\in M_a(\g)$). We now
invoke the weak$^*$-density of $M_a(\g)$ in
${\mathcal P}^*(L_0^\infty(\g)^*)$ to conclude
that ${\mathcal P}^*(L_0^\infty(\g)^*)\odot {\rm\textsf{r}}=0$, where
${\rm\textsf{r}}=\m-\sigma$. That is
${\rm\textsf{r}}\in{\rm ran}({\mathcal P}^*(L_0^\infty(\g)^*))$.$\hfill\blacksquare$\\


In \cite[Page 467]{LOS}, Losert by the use of the C$^*$-algebraic
structure of $M(\g)^*$ proved that $M(\g)^{**}$ has a non-zero
left (weakly) completely continuous element if and only if $\g$ is
compact. Related to this result, we have the following result for
$GL_0(\g)^*$ where our approach in its proof is totally different
from the Losert's result and relies on the theory of generalised
functions.

\begin{theorem}\label{exist}
The following conditions are equivalent.
\newcounter{j1219}
\begin{list}%
{\rm(\roman{j1219})}{\usecounter{j1219}}
\item $\g$ is compact.
\item $GL_0(\g)^*$ has a non-zero left completely continuous element.
\item $GL_0(\g)^*$ has a non-zero left weakly completely continuous element.
\end{list}
\end{theorem}
\noindent Proof. We need only to show that (iii) implies (i).
Indeed, if $\g$ is compact, then $M(\g)^{**}=GL_0(\g)^*$ and the
normalized Haar measure $\m$ on $\g$ is a left (weakly) completely
continuous element of $GL_0(\g)^*$ and (ii)$\Rightarrow$(iii) is
trivial. To this end, suppose that $\m$ is a nonzero left weakly
completely continuous element of $GL_0(\g)^*$. Then the set
$\{\m\odot\delta_x:~x\in\g\}$ is weakly compact and therefore
$\{|\m\odot\delta_x|:~x\in\g\}$ is weakly compact in $GL_0(\g)^*$
by Dieudonne's characterization of weakly compact subsets; see
\cite[Theorem 4.22.1]{DIE}. It follows that ${\cal
E}:=\{|\m|\odot\delta_x:~x\in\g\}$ is weakly compact; this is
because of, $|\m\odot\delta_x|=|\m|\odot\delta_x$ for all
$x\in\g$. Now we apply the Kerin-Smulyan theorem \cite{KS} to
infer that the closed convex hull ${\cal K}$ of ${\cal E}$ is
weakly compact in $GL_0(\g)^*$. On the other hand, it is easy to
see that the map $T_x:{\cal K}\rightarrow {\cal K}$ defined by
$T_x(\n)=\n\odot\delta_x$ is affine for all $x\in\g$. Moreover, we
have
\begin{eqnarray*}
\|\n\odot\delta_x\|&=&\sup{\Big\{}|\big<\n\odot\delta_x,f\big>|:
~f\in GL_0(\g)^*,~\|f\|_\infty\leq 1{\Big\}}\\
&=&\sup{\Big\{}|\big<\n,f\circ\delta_x\big>|:
~f\in GL_0(\g)^*,~\|f\|_\infty\leq 1{\Big\}}\\
&=&\sup{\Big\{}|\big<\n,h\big>|: ~h\in
GL_0(\g)^*,~\|h\|_\infty\leq 1{\Big\}}\\&=&\|\n\|,
\end{eqnarray*}
for all $\n\in GL_0(\g)^*$. It follows that the map $T_x$ is
distal for all $x\in\g$. So there exists a fixed point
${\textsf{q}}\in {\cal K}$ for the maps $T_x$ ($x\in\g$); that is,
${\textsf{q}}\odot\delta_x={\textsf{q}}$ for all $x\in\g$ by the
Ryll-Nardzewski fixed point Theorem; see \cite[Theorem 10.8]{KS}.
In particular, ${\textsf{q}}=\sum_{i=1}^ta_i|\m|\odot\delta_{x_i}$
for some $x_1,...,x_t\in\g$ and $a_1,...,a_t$ with
$\sum_{i=1}^ta_i=1$. Now, if $(K_\alpha)$ denotes the family of
compact subsets of $\g$ ordered by the upward inclusion, then
$(\chi_{{K_\alpha}x^{-1}})$ is a bounded approximate identity for
$GL_0(\g)$ for all $x\in\g$. Thus
\begin{eqnarray*}
\|{\textsf{q}}\odot\delta_x\|
={\bigg\|}\sum_{i=1}^ta_i|\m|\odot\delta_{x_ix}{\bigg\|}
&=&\lim_{\alpha}\sum_{i=1}^ta_i\big<|\m|\odot\delta_{x_ix},\chi_{K_\alpha}\big>\\
&=&\sum_{i=1}^ta_i\lim_{\alpha}\big<|\m|,\chi_{{K_\alpha}{x_ix}^{-1}}\big>
=\|\m\|.
\end{eqnarray*}
Therefore $\|{\textsf{q}}\|=\|\m\|$; since $\m\neq 0$, it follows
that ${\textsf{q}}\neq 0$.

To prove (i), suppose on the contrary that $\g$ is not compact and
that $\overline{{\textsf{q}}}$ is an extension of ${\textsf{q}}$
from $GL_0(\g)$ to a positive functional with the same norm on
$GL(\g)$; see for example \cite{MURP}, Theorem 3.3.8. Then, by the
same manner as in the proof of Proposition \ref{mean}, one can
show that $\overline{{\textsf{q}}}|_{GL_0(\g)}=0$. This implies
that ${\textsf{q}}=0$ a contradiction.
We have now completed the proof of the theorem.$\hfill\blacksquare$\\


We conclude this work by the following result which is of interest in its own right.
In this proposition, the notation
$C_0(\g)^{\perp}$  is used to denote
the set of all $m\in GL_0(\g)^*$ such that $\m|_{C_0(\g)}=0$ and
$${\cal E}_1(\g)={\Big\{}E\in L_0^\infty(\g)^*:~\|E\|=1~{\hbox{and~}}~E~
{\hbox{is~a~right~identity~for~}}L_0^\infty(\g)^*{\Big\}}.$$
It should be noted that $E\in {\cal E}_1(\g)$ if and only if it is a weak$^*$-cluster point
of an approximate identity in $M_a(\g)$ bounded by one; see \cite{LP}.

\begin{proposition}\label{18}
$GL_0(\g)^*$ is
commutative if and only if $\g$ is discrete and abelian.
\end{proposition}
{\noindent Proof}. The necessity of the condition ``$GL_0(\g)^*$ is
commutative" is clear. We prove its sufficiency. To this end, suppose that $GL_0(\g)^*$ is
commutative. That $\g$ is abelian follows trivially.
In order to prove that $\g$ is discrete, we
note that Proposition \ref{05} together with the
fact that the right translations on $GL_0(\g)^*$ are
weak$^*$-continuous implies that
$${\rm ran}({\mathcal P}^*(L_0^\infty(\g)^*))=C_0(\g)^{\perp}.$$
Moreover, by another application
of Proposition \ref{05} one can obtain that ${\mathcal
P}^*(E)-\delta_e\in C_0(\g)^{\perp}$ for
all $E\in {\cal E}_1(\g)$. On the other hand,
from the commutativity of $GL_0(\g)^*$ we get that
$C_0(\g)^{\perp}={\rm lan}({\mathcal P}^*(L_0^\infty(\g)^*))$. We therefore have
$${\mathcal P}^*(E)-\delta_e\in {\rm lan}({\mathcal P}^*(L_0^\infty(\g)^*))\quad
\quad\quad\forall~ E\in {\cal E}_1(\g).$$ It follows that each
elements of ${\mathcal P}^*({\cal E}_1(\g))$ is also a left
identity for ${\mathcal P}^*(L_0^\infty(\g)^*)$. We now invoke
parts (ii) and (iii) of \cite[Theorem 2.11]{LP} to
conclude that $M_a(\g)=M(\g)$. This implies that $\g$ is discrete.$\hfill\blacksquare$\\

\vspace{3mm}

\noindent {\sc Hossein Javanshiri}\\
Department of Mathematics,
Yazd University,
P.O. Box: 89195-741, Yazd, Iran\\
 e-mail: h.javanshiri@yazd.ac.ir\\

\noindent {\sc Rasoul Nasr-Isfahani}\\
Department of Mathematical Sciences,
 Isfahan University of Technology,
  Isfahan 84156-83111, Iran\\
e-mail: isfahani@cc.iut.ac.ir\\

     \end{document}